\documentclass[10pt,a4paper]{article}

\usepackage{mathptmx}       
\usepackage{helvet}         
\usepackage{courier}        
\usepackage{type1cm}        
\usepackage{makeidx}         
\usepackage{graphicx}        
\usepackage{multicol}        
\usepackage{amsmath}

\usepackage{amsfonts,graphicx}
\usepackage{amssymb,hyperref}
\usepackage{color,url}
\usepackage{natbib}

\usepackage[normalem]{ulem}

\makeindex             
\newtheorem{definition}{Definition}

\begin{document}
\begin{center}
\Large{\textbf{Robust Fusion Methods for Big Data}}\\[1.cm]
\end{center}
\normalfont
\centering{Catherine Aaron$^a$, Alejandro Cholaquidis$^b$, Ricardo Fraiman$^b$ and Badih Ghattas$^c$}
\begin{center}
$^a$ Universit\'e Clermont Auvergne, Campus Universitaire des C\'ezeaux, France.\\
$^b$ Universidad de la Rep\'ublica, Facultad de Ciencias, Uruguay.\\
$^c$ Aix Marseille Universit\'e, CNRS, Centrale Marseille, I2M UMR 7373, 13453, Marseille, France.
\end{center}

\abstract{We address one of the important problems in Big Data, namely how to combine estimators from different subsamples by robust fusion procedures, when we are unable to deal with the whole sample.}

\section{Introduction}
Big Data covers a large list of different problems, see for instance \cite{wang2016}, \cite{yu2014}, and references therein. We address one of them, namely how to combine, using robust techniques, estimators obtained from different subsamples in the case where we are unable to deal with the whole sample. In what follows we will refer to this as Robust Fusion Methods (RFM).

To fix ideas, we start by describing one of the simplest problems in this area as a toy example. Suppose we are interested in the median of a huge set of iid  random variables $\{X_1, \ldots, X_n\}$ with common density $f$, and we split the sample into $m$ subsamples of length $l$, $n=ml$. We calculate the median of each subsample and obtain $m$ random variables $Y_1, \ldots, Y_m$. Then we take the median of the set $Y_1, \ldots, Y_m$, i.e. we consider the well known median of medians. It is clear that it does not coincide with the median of the whole original sample $\{X_1, \ldots, X_n\}$, but it will be close. What else could we say about this estimator regarding efficiency and robustness?

The estimator is nothing but the median of $m$ iid random variables but now with a different distribution given by the distribution of the median of $l$ random variables with density $f_X$. Suppose for simplicity that $l=2k+1$. Then, the density of the random variables $Y_i$ is given by 
\begin{equation}\label{densitymedian}
g_Y(y) = \frac{(2k+1)!}{(k!)^2} F_X(t)^k (1-F_X(t))^kf_X(t).  
\end{equation}

In particular, if $f_X$ is uniform on $(0,1)$, (\ref{densitymedian}) is given by	
\begin{equation}\label{densitymedian2}
h_Y(y) = \frac{(2k+1)!}{(k!)^2} t^k (1-t))^k \mathbf 1_{[0,1]}(t),  
\end{equation}
a $Beta(k+1,k+1)$ distribution. 

On the other hand, we have that asymptotically, for a sample of size $n$ the empirical median $\theta=med(X_1, \ldots, X_n)$ behaves as a normal distribution centred at the true median $\theta$ with variance $\frac{1}{4 n f_X(\theta)^2}$, 
while the median of medians behaves asymptotically as a normal distribution centred at $\theta$, the median of the median distribution, and with variance $\frac{1}{4 m g_Y(\theta)^2 }$.
For the uniform case, both are centred at 1/2, $f_X(0.5)=1$ and $g_Y(0.5)=  (1/2)^{2k}(2k+1)!/(k!)^2 $, so we can explicitly calculate the relative efficiency. \\
In section \ref{general} we generalize this idea and study the breakdown point, efficiency, and computational time of the robust fusion method. In section \ref{covop} we tackle, as a particular case, the robust estimation of the covariance operator and show, in section \ref{sim}, the performance throughout a simulation study.


\section{A general setup for RFM.} \label{general}

In this short note we present briefly a general framework for RFM methods for several multivariate and functional data problems. We illustrate our procedure considering only the problem of robust covariance operator estimation, based on a new simple robust estimator. 
Our approach is to consider RFM methods based on data depth functions. The idea is quite simple: given a statistical problem, (such as multivariate location, covariance operators, linear regression, principal components,  among many others), we first split the sample into subsamples.  For each subsample we calculate a robust estimator for the statistical problem considered.  We will use them all to obtain an RFM estimator that is more accurate. More precisely, the RFM estimator is defined as the deepest point (with respect to the appropriate norm associated to the problem) among all the estimators obtained from the subsamples.  Since we need to be able to calculate depths for large sample sizes and eventually high dimensional and infinite dimensional data, we will consider the spatial median corresponding to maximizing  the spatial depth function
\begin{equation}\label{depth}
D(x,P)= 1- \left\|\mathbb E_P \left( \frac{X-x}{\Vert X - x\Vert}\right)\right\|,
\end{equation}
where $P$ is a probability in some Banach space $(E,\|\cdot\|)$ and $x\in E$, introduced by \cite{ch:96},  formulated (in a different way) by \cite{vz:00},  and extended to a very general setup by \cite{chch:14}. We want to address the consistency, efficiency, robustness and computational time properties of the RFM proposals.

To be more precise, the general algorithm is as follows. a) We observe $X_1,\ldots X_n$ iid random elements in a metric space $E$ (for instance $E=\mathbb{R}^d$), b)  we split the sample into subsamples $\{X_1,\ldots X_l\},\{X_{l+1},\ldots X_{2l}\},\ldots,\{X_{(m-1)l+1},\ldots X_{lm}\}$ with $n=ml$, c) we solve our statistical problem on each subsample with a robust procedure (for example, estimate a parameter  $\theta$ on each subsample, obtaining $\hat{\theta}_1,\ldots,\hat{\theta}_m$), d) we take the fusion of the results at each subsample, (for instance $\hat{\hat{\theta}}$ can be the deepest point among $\hat{\theta}_1,\ldots,\hat{\theta}_m$.
	

\subsection{Breakdown point}
\textit{Breakdown point for the RFM.}
Following \cite{do:82} we consider the finite-sample breakdown point. 

\begin{definition} Let $\hat{\theta}_n=\hat{\theta}_n(x)$ be an estimate of $\theta$ defined for samples $\mathbf{x}=\{x_1,\dots,x_n\}$. Let us assume that $\theta$ takes values in $\Theta\subset \mathbb{R}^d$ (it can be $\Theta=\mathbb{R}^d$). Let $\mathcal{X}_p$ be the set of all data sets $\mathbf{y}$ of size $n$ having $n-p$ elements in common with $x$:
	$$\mathcal{X}_p=\{\mathbf{y}:card(\mathbf{y})=n, \ card(\mathbf{x}\cap \mathbf{y})=n-p\},$$
	then $\epsilon_n^*(\hat{\theta}_n,\mathbf{x})=\frac{p^*}{n},$
	where
	$p^*=\max\{p\geq 0: \hat{\theta}_n(\mathbf{y})$ is bounded and also bounded away from $\partial \Theta \ \forall \mathbf{y}\in \mathcal{X}_p\}.$
\end{definition}

Let us consider, for $n=ml$, the random walk $S_n$ with $S_0=0$, and $S_j = B_1+ \ldots + B_j$ for $j=1,\dots,n$, $B_i$ being iid $\text{Bernoulli}(p)$ for $i=1,\dots,n$, where a one represents the presence of an outlier, while a zero represents no presence of an outlier. Then to compute the breakdown point for the median of medians, we need to count how many times the sequence $\{S_l, S_{2l}-S_l, \ldots S_n - S_{n-l}\}$ is larger than $k$ (recall that $l=2k+1$).  Let us define,
$U_{m,n}:= card\{1\leq j \leq m: S_{jl}-S_{(j-1)l}\leq k \} /m,$ 
since the median has breakdown point $0.5$ the fusion will break down if $U_{m,n}$ is greater than $0.5$. 

This will also be true if we take the  median of any robust estimate with breakdown point equal to 0.5 calculated at each subsample.

To have a glance at the breakdown point, we performed 5000 replicates of the vector $S_{30000}$ and calculated the percentage of times the estimator breaks down  for $p=0.45,0.49,0.495$ and $0.499$. The results are in the following table.
\begin{table}[h]
	\begin{center}
		\begin{tabular}{|c|c|c|c||c|}
			\hline
			$m$	& $p=0.45$  & $p=0.49$& $p=0.495$  & $p=0.499$\\
			\hline              
			\hline
			5  & 0  & 0.0020    & 0.0820 & 0.3892  \\
			10 & 0  & 0.0088    & 0.1564 & 0.5352  \\
			30 & 0  & 0.0052    & 0.1426 & 0.5186 \\  
			50 & 0	& 0.0080    & 0.1598 & 0.5412\\    
			100& 0	& 0.0192    & 0.2162 & 0.6084   \\  
			150& 0	& 0.0278    & 0.2728 & 0.6780 \\  
			
			\hline
		\end{tabular}
	\end{center}
	\caption{Percentage of estimator breakdowns for $5000$ replications and different values of $m$ for $n=30000$; $p$ is the proportion of outliers}
\end{table}


%


Since the number $Y$ of outliers in the subsamples of length $l$ follows a Binomial distribution, $ \text{Binom}(l,p)$, as a direct application of Theorem 1 in \cite{sh:13} we can bound the probability, $q=\mathbb{P}(Y>l/2)$, of breakdown.

\subsection{Efficiency of Fusion of $M$-estimators}
In this section we will obtain the asymptotic variance of the RFM method, for the special case of $M$-estimators. Recall that $M$-estimators can be defined (see Section 3.2 in \cite{hu:09}) by the implicit functional equation
$\int \psi(x,T(F))F(dx)=0,$
where $\psi(x;\theta)=(\partial /\partial \theta)\rho(x;\theta)$, for some function $\rho$. The estimator $T_n$ is given by the empirical version of $T$, based on a sample $\mathcal{X}_n=\{X_1,\dots,X_n\}$.
It is well known that $\sqrt{n}(T_n-T(F))$ is asymptotically normal with mean 0, variance $\sigma^2$, and can be calculated in general as the integral of the square of the influence function.
The asymptotic efficiency of $T_n$ is defined as 
${\rm Eff}(T_n)= \frac{\sigma_{ML}^2}{\sigma^2}$, where $\sigma_{ML}^2$ is the asymptotic variance of the maximum likelihood estimator.
Then, the asymptotic variance of a $M$-estimator built from a sample $T_n^1,\dots,T_n^p$ of $p$ $M$-estimators of $T$ can be calculated easily.



\subsection{Computational time}

We want to calculate the computational time of our robust fusion method for a   sample $\mathcal{X}_n=\{X_1,\dots,X_n\}$ iid of $X$, where we have split $\mathcal{X}_n$ into $m$ subsamples of length $l$, then apply a robust estimator to every subsample of length $l$, and fuse them by taking the deepest point among the $m$ subsamples. If we denote by compRE($l$) the computational time required to calculate the robust estimator based on every subsample of length $l$, and compDeph($m$) the computational time required to compute the deepest robust estimator based on the $m$ estimators, then the computational time of our robust fusion method is $m\times$compRE($l$)+compDeph($m$). 

\section{Robust Fusion for covariance operator} \label{covop}

The estimation of the covariance operator of a stochastic process is a very important topic that helps to understand the fluctuations of the random element, as well as to derive the principal functional components from its spectrum. Several robust and non-robust estimators have been proposed, see for instance \cite{chch:14}, and the references therein. In order to perform RFM, we introduce a computationally simple robust estimator to apply to each of the $m$ subsamples, that can be performed using parallel computing. It is based on the notion of impartial trimming applied on the Hilbert--Schmidt space, where covariance operators are defined. The RFM estimator is the deepest point among the $m$ estimators corresponding to each subsample, where the norm in \eqref{depth} is given by \eqref{norm} below. 

\subsection{A resistant estimate of the covariance operator}

Let $E=L^2(T)$, where $T$ is a finite interval in $\mathbb R$, and $X, X_1, \ldots X_n, \ldots$ iid random elements taking values in $(E, \mathcal B(E))$, where $\mathcal B(E)$ stands for the Borel $\sigma$-algebra on $E$. Assume that $\mathbb E(X(t)^2) < \infty$ for all $t\in T$, and $\int_T \int _T \rho^2(s,t) dsdt < \infty$, so the covariance function is well defined and given by
$
\rho(s,t)=\mathbb E((X(t)-\mu(t))(X(s)-\mu(s))), \ \ \mbox{where} \ \mathbb E(X(t))=\mu(t)
$.

For notational simplicity we  assume that $\mu(t)=0,\forall t\in T$. 
Under these conditions, the covariance operator, given by $S_0 (f) = \mathbb E(\langle X,f\rangle X)$,
is diagonalizable, with nonnegative eigenvalues $\lambda_i$ such that $ \sum_i \lambda^2_i < \infty$.  Moreover $S_0$ belongs to the Hilbert--Schmidt space $HS(E,E)$ of linear operators with square norm and inner product given by
\begin{equation}\label{norm}
\Vert S \Vert^2 _{HS} = \sum_{k=1}^{\infty} \Vert S(e_k) \Vert^2 < \infty, \ \ \langle S_1,S_2\rangle_{\mathcal F} = \sum_{k=1}^{\infty} \langle S_1(e_k), S_2(e_k)\rangle,
\end{equation}
respectively, where $\{e_k: k \geq 1\}$ is any orthonormal basis of $E$, and $S,S_1, S_2 \in  HS(E,E)$. 
In particular, $\Vert S_0 \Vert^2 = \sum_{i=1}^{\infty} \lambda_i^2$, where $\lambda_i$ are the eigenvalues of $S_0$.
Given an iid sample $X_1, \ldots,X_n$, define the Hilbert--Schmidt operators of rank one,
$$
W_i: E \rightarrow E, \hspace{2mm}    W_i(f) = \langle X_i, f \rangle X_i(.),  \hspace{2mm}  i=1, \ldots n.
$$
Let $\phi_i = X_i/\Vert X_i\Vert$. Then, $W_i(\phi_i ) = \Vert X_i\Vert^2 \phi_i=: \eta_i \phi_i.$

The standard estimator of $S_0$ is just the average of these operators, i.e. $\hat S_n = \frac{1}{n} \sum_{i=1}^n W_i$, which is a consistent estimator of $S_0$ by the Law of Large Numbers in the space $HS(E,E)$.  Our proposal is to consider an impartial trimmed estimate as a resistant estimator. The notion of impartial trimming was introduced in \cite{go:91}, while the functional data setting was considered in \cite{caf:06}, from were it can be obtained the asymptotic theory for our setting. In order to perform the algorithm we will derive an exact formula for the matrix distances $\Vert W_i - W_j \Vert$, $1 \leq i \leq j \leq n$.

\bf {Lemma 1} \it We have that 
\begin{equation} \label{eqlem1}
d_{ij}^2:=\Vert W_i - W_j \Vert_{HS}^2 =  \Vert X_i\Vert^4 + \Vert X_j\Vert^4 - 2  \langle X_i, X_j\rangle^2 \quad  \text{ for } 1 \leq i \leq j \leq n.
\end{equation}

 \rm



\subsection{The impartial trimmed mean estimator}

Following \cite{caf:06}, we define the impartial trimmed covariance operator estimator, which is calculated by the following algorithm.

Given the sample $X_1(t), \ldots, X_n(t)$  (which we have assumed with mean zero for notational simplicity) and $0<\alpha<1$, we provide a simple algorithm to calculate an approximate impartial trimmed mean estimator of  the covariance operator
of the process $
S_0 : E \rightarrow E,  \hspace{2mm} S_0(f) (t) = \mathbb E(\langle X,f\rangle X(t))$, 
that will be strongly consistent.\\

\textbf{STEP 1:} Calculate $d_{ij}=\Vert W_i - W_j \Vert_{HS}$,  $1 \leq i \leq j \leq n$, using Lemma 1.\\

\textbf{STEP 2:}  Let $r=\lfloor(1-\alpha)n\rfloor+1$.  For each $i=1, \ldots n$, consider the set of indices $I_i \subset \{1, \ldots ,n\}$ corresponding to the $r$ nearest neighbours of $W_i$ among $\{W_1, \ldots W_n\}$, and
the order statistic of the vector $(d_{i1}, \ldots , d_{in})$, $d_i^{(1)}<\ldots <d_i^{(n)}$.\\

\textbf{STEP 3:} Let $\gamma = \text{argmin} \{ d_1^{(r)}, \ldots, d_n^{(r)} \}$.\\

\textbf{STEP 4:} The impartial trimmed mean estimator of $S_0$ is given by 
$\hat S =$ the average of the $m$ nearest neighbours of $W_{\gamma}$ among  $\{ W_1, \ldots, W_n\}$, i.e the average of the observations in $I_\gamma$.

This estimator corresponds to estimating $\rho(s,t)$ by $\hat\rho(s,t) = \frac{1}{r} \sum_{j \in I_\gamma }X_j(s)X_j(t)$. Observe that steps 1 and 2 of the algorithm can be performed using parallel computing.

\subsection{Simulation results for the covariance operator} \label{sim}

Simulations were done using a PC Intel core i7-3770 CPU, 8GO of RAM using 64 bit version of Win10, and R software ver. 3.3.0.

We vary  the sample size $n$ within the set $\{0.1e6,1e6,5e6,10e6\}$ and the number of subsamples $m\in \{100,500,1000,10000 \}$. The proportion of outliers was fixed to $p=13\%$ and $p=15\%$. We replicate each simulation case $K=5$ times and report a mean average of the results over these replicates.

 We report the average time in seconds necessary for both a global estimate (\textit{time0}, over the whole sample), and \textit{time1} the estimate obtained by fusion (including computing the estimates over subsamples and aggregating them by fusion). \\
We compare the classical estimator (\textsl{Cov}), the mean of the classical estimators obtained from the subsamples (\textsl{AvCov}), the Fusion estimate of the classical estimator (\textit{Cov.RFM}), the global robust estimate (\textit{CovRob}), the robust fusion estimate \textit{RFM}, and the average of the robust estimates from the subsamples \textit{AvRob}.

To generate the data, we have used a simplified version of  the simulation model used in \cite{Kraus}:
$$ X(t) = \mu(t) + \sqrt 2 \sum_{k=1}^{10} \lambda_k a_k \sin(2 \pi k t) + \sqrt 2 \sum_{k=1}^{10} \nu_k b_k \cos(2 \pi k t)$$
where $\nu_k = \left(\frac{1}{3}\right) ^k, \lambda_k = k^{-3}$, and $a_k$ and $b_k$ are random standard Gaussian independent observations. The central observations were generated using $\mu(t) = 0$ whereas for the outliers we took $\mu(t) =2- 8\sin(\pi t)$. For $t$ we used an equally spaced grid of $T=20$ points in $[0,1]$.\\
The covariance operator of this process was computed for the comparisons: 

$ Cov(s,t)= \sum_{k=1}^{10} A_k(s)A_k(t) + B_k(s)B_k(t) $, 
where $A_k(t) = \sqrt 2 \lambda_k \sin(2 \pi k t)$ and $ A_k(t) = \sqrt 2 \nu_k \cos(2 \pi k t)$.

The results are shown in the following two tables for two proportions of outliers, $p=0.15$ and $p=0.2$.
\begin{table}[h!]
	\caption{Covariance operator estimator in presence of outliers. Using the classical and robust estimators over the entire sample, and aggregating by average or fusion of M subsamples estimates. p=0.15, T=20\label{signif}} 
	\begin{center}
		\begin{tabular}{rrrrrrrrrr}
			\hline\hline
			\multicolumn{1}{c}{n}&\multicolumn{1}{c}{M}&\multicolumn{1}{c}{time0}&\multicolumn{1}{c}{time1}&\multicolumn{1}{c}{Cov}&\multicolumn{1}{c}{AvCov}&\multicolumn{1}{c}{Cov.RFM}&\multicolumn{1}{c}{CovRob}&\multicolumn{1}{c}{AvRob}&\multicolumn{1}{c}{RFM}\tabularnewline
			\hline
			$0.05$&$  20$&$ 553$&$18.20$&$24.3$&$24.3$&$24.7$&$5.16$&$5.21$&$5.52$\tabularnewline
			$0.05$&$  50$&$ 543$&$ 7.81$&$24.3$&$24.3$&$24.9$&$5.20$&$5.24$&$5.60$\tabularnewline
			$0.05$&$ 100$&$ 528$&$ 4.79$&$24.3$&$24.3$&$25.2$&$5.20$&$5.17$&$5.58$\tabularnewline
			$0.05$&$1000$&$ 459$&$19.40$&$24.3$&$24.3$&$27.0$&$5.13$&$5.54$&$6.58$\tabularnewline
			$0.10$&$  20$&$2300$&$69.00$&$24.2$&$24.2$&$24.4$&$5.14$&$5.22$&$5.43$\tabularnewline
			$0.10$&$  50$&$2300$&$28.10$&$24.2$&$24.2$&$24.6$&$5.04$&$5.09$&$5.13$\tabularnewline
			$0.10$&$ 100$&$2290$&$15.20$&$24.2$&$24.2$&$25.0$&$5.06$&$5.15$&$5.43$\tabularnewline
			$0.10$&$1000$&$1850$&$21.60$&$24.3$&$24.2$&$26.1$&$5.21$&$5.35$&$6.13$\tabularnewline
			\hline
	\end{tabular}\end{center}
	
\end{table}

\begin{table}[h!]
	\caption{Covariance operator estimator in the presence of outliers. Using classical and robust estimators over the entire sample, and aggregating by average or fusion of $m$ subsamples estimates. $p=0.2$, T=20\label{signif}} 
	\begin{center}
		\begin{tabular}{rrrrrrrrrr}
			\hline\hline
			\multicolumn{1}{c}{n}&\multicolumn{1}{c}{$m$}&\multicolumn{1}{c}{time0}&\multicolumn{1}{c}{time1}&\multicolumn{1}{c}{Cov}&\multicolumn{1}{c}{AvCov}&\multicolumn{1}{c}{Cov.RFM}&\multicolumn{1}{c}{CovRob}&\multicolumn{1}{c}{AvRob}&\multicolumn{1}{c}{RFM}\tabularnewline
			\hline
			$0.05$&$  20$&$ 572$&$17.90$&$30.5$&$30.5$&$30.9$&$0.879$&$ 3.96$&$1.45$\tabularnewline
			$0.05$&$  50$&$ 649$&$ 7.88$&$30.5$&$30.5$&$31.3$&$0.876$&$ 7.34$&$2.10$\tabularnewline
			$0.05$&$ 100$&$ 633$&$ 4.61$&$30.5$&$30.5$&$31.6$&$0.839$&$ 8.86$&$2.43$\tabularnewline
			$0.05$&$1000$&$ 478$&$19.50$&$30.5$&$30.5$&$32.3$&$0.864$&$13.10$&$7.08$\tabularnewline
			$0.10$&$  20$&$1970$&$69.10$&$30.4$&$30.4$&$30.6$&$0.914$&$ 3.83$&$1.36$\tabularnewline
			$0.10$&$  50$&$2030$&$28.10$&$30.4$&$30.4$&$31.1$&$0.921$&$ 4.32$&$1.55$\tabularnewline
			$0.10$&$ 100$&$2020$&$15.10$&$30.4$&$30.4$&$31.3$&$0.840$&$ 8.44$&$2.35$\tabularnewline
			$0.10$&$1000$&$1840$&$21.60$&$30.4$&$30.4$&$32.9$&$0.961$&$12.10$&$5.20$\tabularnewline
			\hline
	\end{tabular}\end{center}
	
\end{table}

If the proportion of outliers is moderate $p=15\%$, the average of the robust estimators still behaves well, better than RFM, but if we increase the proportion of outliers to $p=0.2$, RFM clearly outperforms all the other estimators. 



\section*{Acknowledgement} 
We thanks an anonymous referee for helpful suggestions on a first version.

\end{document}